\documentclass[11pt, twoside, eqno]{article}
\usepackage[bookmarksnumbered, colorlinks, plainpages]{hyperref}
\usepackage{amssymb}
\usepackage{mathrsfs}
\usepackage{amsmath}
\usepackage{enumitem}
\usepackage{amsfonts}
\usepackage{colortbl,dcolumn}
\usepackage{color, soul}
\usepackage{amsthm}
\hypersetup{colorlinks=true,linkcolor=red, anchorcolor=green, citecolor=cyan, urlcolor=red, filecolor=magenta, pdftoolbar=true}

\allowdisplaybreaks \numberwithin{equation}{section}

\newtheorem{rem}{\noindent\rm\bf Remark}[section]
\newtheorem{cor}{\noindent\rm\bf Corollary}[section]
\newtheorem{thm}{\noindent\rm\bf Theorem}[section]
\newtheorem{lem}{\noindent\rm\bf Lemma}[section]

\newtheorem*{thmA}{Theorem A}

\begin{document}

\baselineskip=15pt
\renewcommand{\arraystretch}{2}
\arraycolsep=1pt
\linespread{1.1}

\title
{\bf \Large Strong unique continuation for the $\bar{\partial}$-operator with $L^2$ potentials
\footnotetext{\hspace{-0.25cm}
\noindent{{*Corresponding Author.
          }\\
\date{}}}}

\author
{\vspace{-0.5cm}\bf
 Xiaomin Tang, Shaoyong He*
\\
\vspace{-0.5cm}\small\it  Department of Mathematics, Huzhou  Normal University,  Huzhou,  313000, P. R. China;\\
\vspace{-0.5cm}\small\it E-mails: ~txm@huznu.edu.cn, hsyongmath@sina.com \\
}

\maketitle

\begin{center}
\begin{minipage}{15.2cm}
\small{\bf Abstract}\quad
In this paper, we prove strong unique continuation for solutions of the differential
inequality $|\bar\partial u|\leq V|u|$ in the complex plane.
If $V\in L_{\mathrm{loc}}^2(\mathbb C)$
and $1<p<2$, then any solution $u\in W_{\mathrm{loc}}^{1,p}(\mathbb C)$ that vanishes to infinite order at a point in the $L^p$ sense must vanish identically.
This result weakens the hypotheses of an earlier unique continuation theorem (Trans. Amer. Math. Soc. 318 (1990), 275--300).
The proof uses a Carleman inequality derived from a sharp $L^p$-$L^q$ estimate for the $\bar{\partial}$-operator via Calder\'{o}n-Zygmund theory.
We also establish the sharpness of the
$L^2$ integrability condition by constructing
a nontrivial smooth counterexample through
a gluing procedure involving Cauchy kernels,
with a potential belonging to $L^p(\mathbb C)$
for every $0<p<2$.

\medskip

\noindent{\bf Keywords}\quad $\bar{\partial}$-operator; unique continuation; Calder\'{o}n-Zygmund theory; Carleman inequality

\noindent{\bf MR(2020) Subject Classification}\quad 32W05, 35A02

\end{minipage}
\end{center}

 \section{Introduction}
Let \(u=(u_1,\ldots,u_N)\) be a vector-valued function defined on a domain \(\Omega\subset\mathbb{C}\), taking values in \(\mathbb{C}^N\).
 We say that a differential equation or inequality enjoys the strong unique continuation property if any solution that vanishes to infinite order at some point $z_0\in\Omega$ must vanish identically throughout $\Omega$.
 A function $u\in L_{\mathrm{loc}}^{p}(\Omega)$ is said to vanish to infinite order (or simply be flat) at $z_0\in\Omega$ if,
 for every $m\geq0$,
\begin{align}\label{1.1}
\lim_{r\rightarrow0}r^{-m}\int_{|z-z_0|<r}|u(z)|^{p}\,dv_z=0,
\end{align}
where $dv_z$ denotes the Lebesgue measure on $\mathbb{C}$.
A differential equation or inequality is said to possess the weak unique continuation property if every solution that vanishes on a nonempty open subset is identically zero.

In the course of studying the boundary regularity and uniqueness of holomorphic maps sending a hypersurface into a Levi-flat hypersurface, Bell and Lempert \cite{BL} proved a strong unique continuation theorem for solutions of the differential inequality
 $|\bar{\partial}u |\leq C|u|$ with $C$ a constant.
Pan and Zhang \cite{PZ} investigated the strong unique continuation property for the more general differential inequality
\begin{align}\label{1.2}
|\bar{\partial}u|\leq V|u|,
\end{align}
where the potential $V\in L_{\mathrm{loc}}^{2}(\Omega) $.
We write $W^{k,p}_{\mathrm{loc}}(\Omega)$ for the space of functions $u\in L^p_{\mathrm{loc}}(\Omega)$ whose weak derivatives of all orders up to $k
 $ belong to $L^p_{\mathrm{loc}}(\Omega)$.
 Pan and Zhang proved the following strong unique continuation theorem.
\begin{thmA}
Let $\Omega\subset\mathbb C$ be a domain, and suppose that
$u\in W^{1,2}_{\mathrm{loc}}(\Omega)$ satisfies \eqref{1.2} with
$V\in L^2_{\mathrm{loc}}(\Omega)$.
If $u$ vanishes to infinite order at some point $z_0\in\Omega$
in the $L^2$ sense, then $u\equiv 0$ in $\Omega$.
\end{thmA}
The proof of Theorem~A is based on a pointwise identity and a weighted inequality for the Riesz potential.
In the present paper, we develop a simpler approach using Calder\'{o}n-Zygmund theory, which also enables us to strengthen and extend Theorem~A.
\begin{thm}\label{thm1}
Let $u\in W^{1,p}_{\mathrm{loc}}(\Omega)$ with $1<p<2$. Suppose that
\[|\bar{\partial}u|\leq V|u|\]
almost everywhere in $\Omega$, where $V\in L^2_{\mathrm{loc}}(\Omega)$.
If $u$ is flat at a point $z_0\in\Omega$ in the $L^p$ sense, then $u$ vanishes identically in $\Omega$.
\end{thm}
Observe that
$$
W^{1,2}_{\mathrm{loc}}(\Omega)
\subset W^{1,p}_{\mathrm{loc}}(\Omega),
\qquad 1<p<2.
$$
Moreover, H\"older's inequality shows that $L^2$-flatness at a point
implies $L^p$-flatness for every $1<p<2$. Thus, the assumptions of
Theorem~A imply those of Theorem~\ref{thm1}, which extends
Theorem~A to solutions with lower Sobolev regularity.

The strong unique continuation property for second-order elliptic
equations has been extensively investigated.
In particular, Jerison and Kenig~\cite{JK1985} established
strong unique continuation for solutions
$u\in W_{\mathrm{loc}}^{2,\frac{2d}{d+2}}(\Omega)$
of the differential inequality
$$
|\Delta u|\leq V|u|,
$$
where $\Omega\subset\mathbb R^d$ is an open connected set
with $d\ge3$ and $V\in L_{\mathrm{loc}}^{d/2}(\Omega)$.
Sogge~\cite{S1990} extended this theory to second-order
elliptic operators $P(x,D)$ with smooth coefficients and a
positive principal symbol in dimensions $d\ge3$.
He established strong unique continuation for the differential inequality
$$
|P(x,D)u|\leq |Vu|+C|x|^{-1+\sigma}|\nabla u|,
$$
where $V\in L_{\mathrm{loc}}^{n/2}(\Omega)$ and $\sigma>0$,
for solutions satisfying $u\in W_{\mathrm{loc}}^{1,q}(\Omega)$
and $P(x,D)u\in L_{\mathrm{loc}}^q(\Omega)$, with $q=2d/(d+2)$.

On the other hand, Wolff~\cite{wolff1990} investigated
strong unique continuation for the differential inequality
\begin{equation}\label{Delta1}
|\Delta u|\leq V|\nabla u|
\end{equation}
in dimensions $d\ge3$, under suitable integrability
assumptions on $V$.
Chanillo and Sawyer~\cite{CS}
established strong unique continuation for \eqref{Delta1}
on a domain $\Omega\subset\mathbb R^2$ with
$V\in L_{\mathrm{loc}}^2(\Omega)$.
More precisely, they proved that any solution
$u\in W_{\mathrm{loc}}^{2,2}(\Omega;\mathbb R^N)$
that vanishes to infinite order at a point $z_0\in\Omega$
must vanish identically.
In fact, by applying Theorem~\ref{thm1} to \(\partial u\), we can relax the regularity
assumption on \(u\) to \(W_{\mathrm{loc}}^{2,p}(\Omega)\) with \(1<p<2\).
This reduction is possible because \(\Delta = 4\bar\partial\partial\) and, since
\(u\) is real-valued, we have the pointwise identity
\(|\partial u(z)| = \frac{1}{2}|\nabla u(x,y)|\).
Thus the original second-order inequality becomes a first-order inequality for
\(\partial u\), which falls under the scope of Theorem~\ref{thm1}.
Consequently, the Chanillo-Sawyer result is recovered as a special case, while
our theorem significantly lowers the required regularity.

In addition to the local vector-valued setting considered in
Theorem~\ref{thm1}, unique continuation for \eqref{1.2}
has also been investigated in the global setting on
$\mathbb C$, particularly for scalar-valued solutions.
We refer the reader to \cite{kenig,PZ2,seo} for related results
with $L^2(\mathbb C)$ potentials.

We next examine the sharpness of the $L^2$ integrability condition on the potential in Theorem~\ref{thm1}.
Mandache~\cite{MAND} constructed a counterexample showing that weak unique continuation may fail when $V\in L^p$ for $0<p<2$. Inspired by Wolff's construction~\cite{wolff1992}, we develop an alternative approach based on Cauchy kernels and a gluing procedure.
This construction yields a nontrivial smooth solution that vanishes to infinite order at the origin, with a potential belonging to $L^p(\mathbb C)$ for all $0<p<2$ simultaneously. More precisely, we obtain the following theorem.

\begin{thm}\label{thm2}

There exist a nontrivial function
$u\in C^\infty(\mathbb C;\mathbb C^N)$
and a nonnegative potential
$
V\in\displaystyle\bigcap_{0<p<2}L^p(\mathbb C)
$
such that
\[
|\bar\partial u|\leq V|u|.
\]
Moreover, $u$ vanishes to infinite order
at the origin.

\end{thm}

\begin{rem}
Unlike the previously known compactly supported
counterexamples, our construction provides a
different approach based on gluing holomorphic
functions across shrinking annuli and
regularizing isolated poles.
The main feature of this method is that it
allows separate control of the gluing errors
and the singularities of the potential,
while ensuring infinite-order vanishing
at the origin. This provides an alternative
perspective on the sharpness of the
critical integrability exponent $2$. We also note that it suffices to prove
Theorem~\ref{thm2} for $N=1$.
Indeed, suppose that a nontrivial
scalar-valued function $u\in C^\infty(\mathbb C)$
vanishes to infinite order at the origin
and satisfies $|\bar\partial u|\leq V|u|$
for some $V\in L^p(\mathbb C)$, $0<p<2$. For any $N\geq1$, define
\[
\widetilde u=(u,0,\ldots,0)
\in C^\infty(\mathbb C;\mathbb C^N).
\]
Then
\[
|\bar\partial\widetilde u|
=|\bar\partial u|
\leq V|u|
=V|\widetilde u|.
\]
Since $\widetilde u$ also vanishes to infinite
order at the origin, the scalar counterexample
extends to every $N\geq1$.
\end{rem}

 \section{Proof of Theorem~\ref{thm1} }
 We begin by establishing a key $L^p$-$L^q$ estimate for the Cauchy-Riemann operator. This estimate is the sharp Sobolev inequality in $\mathbb{R}^2$ for the complex derivative.
 \begin{lem}\label{lem2.1}
Let $u \in C_0^\infty(\mathbb{C})$. If $1<p<2<q<\infty$ and $1/p - 1/q = 1/2$, then
\begin{align}\label{2.1}
\|u\|_{L^q} \leq C \|\bar{\partial} u\|_{L^p}.
\end{align}
\end{lem}
\begin{proof}
Let $h = \bar{\partial} u$. We may  consider the operator $\partial_x + i\partial_y$. Taking the Fourier transform gives
\[
\widehat{(\partial_x + i\partial_y)u}(\xi,\eta) = 2\pi i(\xi + i\eta) \hat{u}(\xi,\eta).
\]
After dropping the constant $2\pi i$, we have
\[
(\xi + i\eta) \hat{u}(\xi,\eta) = \hat{h}(\xi,\eta).
\]
Thus, proving \eqref{2.1} is equivalent to establishing the boundedness of the multiplier operator
$$
T: L^p \to L^q, \quad \widehat{T h}(\xi,\eta) = \frac{1}{\xi + i\eta} \hat{h}(\xi,\eta).
$$
The difficulty is that the multiplier $m(\xi,\eta) = \frac{1}{\xi + i\eta}$ has a non-integrable singularity at the origin, since $|m|^2$ is not integrable near $(0,0)$.


To circumvent this singularity, we employ the Littlewood-Paley theorem. Let $\psi \in C_0^\infty(\mathbb{R}^2)$ be supported in  $\{1/2 < |(\xi,\eta )| < 2\}$ satisfying
\begin{align}\label{sum}
\sum_{k \in \mathbb{Z}} \psi(2^k \xi, 2^k \eta )^2 = 1, \qquad \xi \neq 0.
\end{align}
Define  $\psi_k(\xi,\eta) = \psi(2^k \xi, 2^k \eta )$, and let  $S_k$  denote the associated Littlewood-Paley projection operators
$$
\widehat{S_k h}(\xi,\eta) = \psi_k(\xi,\eta) \hat{h}(\xi,\eta).
$$
The Littlewood-Paley theorem tells us  that for $1<p<\infty$,
$$
\|h\|_{L^p} \approx \left\| \left( \sum_{k \in \mathbb{Z}} |S_k h|^2 \right)^{1/2} \right\|_{L^p}.
$$

Next, we define the localized multiplier operators $T_k$ by
\[
\widehat{T_k h}(\xi,\eta) = m(\xi,\eta) \psi_k(\xi,\eta) \hat{h}(\xi,\eta).
\]
Equivalently,
\[
T_k h = K_k * h,
\]
where
$
K_k = \mathcal{F}^{-1}\!\left( \frac{ \psi_k(\xi,\eta)}{\xi + i\eta} \right).
$
We claim that
\begin{align}\label{2.4}
\|T_k h\|_{L^q} \leq C \|h\|_{L^p}
\end{align}
uniformly in $k$. Indeed, by Young's inequality,
\[
\|T_k h\|_{L^q} \leq \|K_k\|_{L^2} \|h\|_{L^p},
\]
where $1 + \frac 1 q = \frac 1 2 + \frac 1 p$. By Plancherel's theorem,
\begin{align*}
\|K_k\|_{L^2}^2&= \int_{\mathbb{R}^2} \left| \frac{ \psi_k(\xi,\eta)}{\xi + i\eta} \right|^2 \, d\xi d\eta \\
&\leq \int_{\mathbb{R}^2} \frac{|\psi_k(\xi,\eta)|^2}{\xi^2 + \eta^2} \, d\xi d\eta \\
&= \int_{A_k} \frac{|\psi(2^k \xi, 2^k \eta)|^2}{\xi^2 + \eta^2} \, d\xi d\eta,
\end{align*}
where $A_k = \{2^{-k-1} \leq |(\xi,\eta)| \leq 2^{-k+1}\}$. Making the change of variables $(\xi,\eta) = 2^{-k}(\zeta,\omega)$, we obtain
\[
\|K_k\|_{L^2}^2 \leq \int_{A_0} \frac{|\psi(\zeta,\omega)|^2}{\zeta^2 + \omega^2} \, d\zeta d\omega \leq C,
\]
independent of $k$. This proves \eqref{2.4}.

Using duality and  \eqref{sum}, we obtain
\begin{align*}
\|Th\|_{L^q} &= \sup_{\|g\|_{L^{q'}}\leq 1 }\left| \int_{\mathbb{R}^2} (Th)(x) \overline{g(x)} \, dx \right|\\
&=\int_{\mathbb{R}^2} \sum_{k} T_k(S_k h)(x) \overline{\sum_j (S_j S_j g)(x)}dx\\
 &= \sum_{k, j} \int_{\mathbb{R}^2} T_k(S_k h)(x) \overline{(S_j S_j g)(x)} dx.
\end{align*}
Most of the terms in the double sum vanish by Parseval's identity, since
$$\int_{\mathbb{R}^2} T_k(S_k h)(x) \overline{(S_j S_j g)(x)}dx=\int_{\mathbb{R}^2} \widehat{T_k(S_k h)}(\xi,\eta)\overline{\widehat{(S_j S_j g)}}(\xi,\eta),$$
and this integral is zero unless the Fourier supports of $\widehat{T_k(S_k h)}$ and  $\widehat{S_j S_jg}$ intersect. Therefore,
\begin{align*}
    &\sum_{k, j} \int_{\mathbb{R}^2} T_k(S_k h)(x) \overline{(S_j S_j g)(x)} dx   \\
    &= \sum_{k \in \mathbb{Z}} \sum_{|j-k| \le 2} \int T_k(S_k h)(x) \overline{S_j S_j g(x)} \, dx   \\
    &= \sum_{k \in \mathbb{Z}} \sum_{\ell=-2}^{2} \int T_k(S_k h)(x) \overline{S_{k+\ell}S_{k+\ell} g(x)} \, dx  \\
    &= \sum_{\ell=-2}^{2} \sum_{k \in \mathbb{Z}} \int T_k(S_k h)(x) \overline{S_{k+\ell}S_{k+\ell} g(x)} \, dx.
\end{align*}
For each fixed $\ell$, H\"{o}lder's inequality, the estimate \eqref{2.4}, the Minkowski's inequality, and the Littlewood-Paley theorem yield
\begin{align*}
    & \left| \sum_{k \in \mathbb{Z}} \int T_k(S_k h)(x) \overline{S_{k+\ell} S_{k+\ell}g(x)} \, dx \right| \\
    & \le \sum_k \left( \int |T_k(S_k h)(x)|^q dx \right)^{\frac{1}{q}} \left( \int |S_{k+\ell} g(x)|^{q'}dx \right)^{\frac{1}{q'}} \\
    & \le C \sum_k \|S_k h\|_{L^p} \|S_{k+\ell} g\|_{L^{q'}} \\
    & \le C\left( \sum_k \|S_k h\|_{L^p}^2 \right)^{\frac{1}{2}} \left( \sum_k \|S_{k+\ell} g\|_{L^{q'}}^2 \right)^{\frac{1}{2}}\\
    & \le C\| h\|_{L^p} \| g\|_{L^{q'}}.
\end{align*}
This  completes the proof of Lemma \ref{lem2.1}.

\end{proof}
For the vector-valued case $N>1$, we provide a new, shorter, and more general proof of Lemma \ref{lem2.1} using the classical Calder\'{o}n-Zygmund theory of singular integrals.

\begin{lem}\label{lem2.2}
Let $u=(u_1,\dots,u_N)\in C_0^\infty(\mathbb C;\mathbb C^N)$. If $1<p<2<q<\infty$ and
$
\frac{1}{p}-\frac{1}{q}=\frac{1}{2},
$
then there exists a constant $C>0$, independent of $u$, such that
\begin{equation*}
\|u\|_{L^q} \le C \|\bar{\partial} u\|_{L^p},
\end{equation*}
where
$$
\|u\|_{L^q}=\big\|(\sum_{k=1}^N |u_k|^2)^{1/2}\big\|_{L^q},
$$
and similarly for $\bar{\partial}u$.
\end{lem}
\begin{proof}
We recall the Beurling transform $\mathcal B$, defined on scalar functions $f\in C_0^\infty(\mathbb C)$ by the principal-value singular integral
\[
\mathcal B f(z)=\mathrm{p.v.}\int_{\mathbb C}\frac{f(\zeta)}{(z-\zeta)^2}\,d\zeta\,d\bar\zeta.
\]
The kernel $K(z)=1/( z^2)$ is a standard Calder\'{o}n-Zygmund kernel on $\mathbb R^2$: it satisfies $|K(z)|\le C|z|^{-2}$ and $|\nabla K(z)|\le C|z|^{-3}$.
Hence, by the classical Calder\'{o}n-Zygmund theorem,
\begin{equation}\label{vector}
\|\mathcal B u\|_{L^r}=\big\|(\sum_{k=1}^N |\mathcal B u_k|^2)^{1/2}\big \|_{L^r} \le C \big\|(\sum_{k=1}^N |u_k|^2)^{1/2}\big\|_{L^r}=C\|u\|_{L^r}, 1<r<\infty.
\end{equation}
As is well-known, the standard identity $\partial u_k=\mathcal B(\bar{\partial}u_k)$ holds for compactly supported smooth functions.
Therefore, we have $\partial u=\mathcal B(\bar{\partial}u)$. Applying \eqref{vector} with $r=p$ gives
\[
\|\partial u\|_{L^p} \le C\|\bar{\partial}u\|_{L^p}.
\]
For each component $u_k$, using $\partial_x=\partial+\bar\partial$ and $\partial_y=-i(\partial-\bar\partial)$, we obtain
\[
|\partial_x u_k|^2+|\partial_y u_k|^2 = 2\bigl(|\partial u_k|^2+|\bar\partial u_k|^2\bigr).
\]
Summing over $k=1,\dots,N$ and taking square roots yields the pointwise identity
\[
|\nabla u| = \sqrt{2}\left(\sum_{k=1}^N |\partial u_k|^2+\sum_{k=1}^N |\bar\partial u_k|^2\right)^{1/2},
\]
where $|\nabla u|=(\sum_k(|\partial_x u_k|^2+|\partial_y u_k|^2))^{1/2}$. Consequently,
\begin{align}\label{nabla}
\|\nabla u\|_{L^p} \le C \big(\|\bar\partial u\|_{L^p}+\|\partial u\|_{L^p}\big)\leq C\|\bar{\partial}u\|_{L^p}.
\end{align}

It remains to apply the Sobolev embedding theorem to vector-valued functions. Define the scalar function $v(z)=|u(z)|=(\sum_{k=1}^N |u_k(z)|^2)^{1/2}$. Since $u\in C_0^\infty$, $v$ is Lipschitz and compactly supported, hence $v\in W^{1,p}(\mathbb R^2)$. Moreover, we have the pointwise estimate
\[
|\nabla v| \le |\nabla u| \quad \text{a.e.}
\]
Indeed, at points where $u\neq 0$,
\[
\partial_j v = \frac{1}{v}\sum_{k=1}^N \operatorname{Re}(\overline{u_k}\,\partial_j u_k),
\]
so by Cauchy-Schwarz, $|\partial_j v| \le (\sum_k|\partial_j u_k|^2)^{1/2}$; summing over $j=x,y$ gives the desired inequality. On the zero set of $u$,
we  note that the set
\[
E=\{z\in\mathbb C: u(z)=0,\ \nabla u(z)\neq 0\}
\]
has Lebesgue measure zero.
Applying the classical Sobolev embedding theorem in dimension two to the scalar function $v$, with $1<p<2$ and $\frac 1 p-\frac 1 q=\frac 1 2$, yields
\[
\|v\|_{L^q} \le C \|\nabla v\|_{L^p}.
\]
Since $\|u\|_{L^q}=\big\|(\sum_{k=1}^N |u_k|^2)^{1/2}\big\|_{L^q}=\|v\|_{L^q}$, we obtain
\begin{align}\label{sobolev}
\|u\|_{L^q(\ell^2)} \le C \|\nabla u\|_{L^p(\ell^2)}.
\end{align}
Combining (\ref{nabla}) and (\ref{sobolev}), we conclude
\[
\|u\|_{L^q} \le C \|\bar\partial u\|_{L^p},
\]
which is the desired estimate.

\end{proof}
\begin{rem}
The  proof of Lemma~\ref{lem2.1}  relied on Littlewood-Paley theory, which gives the estimate componentwise.
However, that approach does not extend directly to vector-valued functions: the componentwise inequalities $\|u_k\|_{L^q}\le C\|\bar\partial u_k\|_{L^p}$ cannot be combined to yield the desired mixed-norm estimate $\|u\|_{L^q}\le C\|\bar\partial u\|_{L^p}$.
In contrast, the present proof uses the Beurling transform as a Calder\'{o}n--Zygmund operator, which is  bounded on $L^r$, $1<r<\infty$. Together with the pointwise equivalence $|\nabla u|\approx|\partial u|+|\bar\partial u|$ and the Sobolev embedding applied to the scalar function $|u|$, it yields the vector-valued estimate directly, without any componentwise summation. This approach is not only shorter but also genuinely covers the case $N>1$.

\end{rem}

The following Carleman inequality is an immediate consequence of Lemma \ref{lem2.2} and will be the key ingredient in the proof of the unique continuation theorem.
\begin{cor}\label{cor2.2}
Let $1<p<2<q<\infty$ with $\frac 1 p - \frac 1 q = \frac 1 2$.
For $u \in C_0^\infty(\mathbb{C} \setminus \{0\}, \mathbb{C}^{N} )$ and integer $\tau > 0$,
\begin{align}\label{2.21}
\left\| |z|^{-t} u \right\|_{L^q} \leq C \left\| |z|^{-t} \bar{\partial} u \right\|_{L^p}.
\end{align}
\end{cor}
For $z \in \mathbb{C} \setminus \{0\}$, note that
$$
\bar{\partial}(z^{-t} u) =z^{-t}\bar{\partial} u+u \bar{\partial}(z^{-t}) =z^{-t} \bar{\partial} u.
$$
Thus, setting $f = z^{-t} u$, inequality \eqref{2.21} is equivalent to
$$
\|f\|_{L^q} \leq C \|\bar{\partial} f\|_{L^p}.
$$
But this is precisely Lemma \ref{lem2.2} applied to $f$, since $f \in C_0^\infty(\mathbb{C} \setminus \{0\}, \mathbb{C}^{N} )$.

  We now turn to the
\begin{proof}[Proof of Theorem \ref{thm1}.] Without loss of generality, we assume $z_0 = 0$.
Let $\varphi \in C_0^\infty(\mathbb{D})$ satisfy
\[
0\leq \varphi\leq 1,\qquad
\varphi\equiv 1\quad \text{on } B(0,\tfrac1 2),
\]
and set
$
f=\varphi u.
$
Since $u\in W^{1,p}_{\mathrm{loc}}$, we have
$f\in W^{1,p}(\mathbb{D})$.

Let $\{\Psi_j\}_{j=1}^{\infty}\subset C^\infty(\mathbb{D})$ be a sequence satisfying
\[
0\leq \Psi_j\leq \Psi_{j+1}\leq 1,
\]
where
\[
\Psi_j=0\quad\text{for } |z|<\frac1j,\qquad
\Psi_j=1\quad\text{for } |z|>\frac2j,
\]
and
$
|\nabla \Psi_j|\leq Cj .
$
Choose $\rho>0$ with $\rho<\frac12$ such that
\[
C\|V\|_{L^2(B(0,\rho))}<\frac1 2 .
\]
We shall prove that $u\equiv0$ in $B(0,\rho)$.

Let
$
f_j=\Psi_j f.
$
By a standard approximation, Corollary ~\ref{cor2.2} holds for $f_j\in W^{1,p}(\mathbb{D} \backslash \{0\}).$
Therefore,
$$
\||z|^{-\tau}f_j\|_{L^q(\mathbb C)}
\leq
C\||z|^{-\tau}\bar\partial f_j\|_{L^p(\mathbb C)}.
$$
Since
$
\bar\partial f_j
=
(\bar\partial\Psi_j)f+\Psi_j\bar\partial f,
$
we obtain
\begin{align*}
\||z|^{-\tau}\Psi_j u\|_{L^q(B(0,\rho))}&=\||z|^{-\tau}\Psi_j \varphi u\|_{L^q(B(0,\rho))}\\
&=\||z|^{-\tau}f_j\|_{L^q(\mathbb C)}\\
&\leq
C\||z|^{-\tau}(\bar\partial\Psi_j)f\|_{L^p(\mathbb C)}
\\
&\quad
+C\||z|^{-\tau}\Psi_j\bar\partial f\|_{L^p(|z|<\rho)}+\||z|^{-\tau}\Psi_j\bar\partial f\|_{L^p(|z|>\rho)}.
\end{align*}

Because $\varphi\equiv1$ on $B(0,\rho)$, we have $f=u$ in $B(0,\rho)$ and hence
$
\bar\partial f=\bar\partial u$
in $B(0,\rho).$
Using the assumption
$
|\bar\partial u|\leq V|u|,
$
we get
\begin{align*}
\||z|^{-\tau}\Psi_j\bar\partial f\|_{L^p(|z|<\rho)}
&\leq
\||z|^{-\tau}\Psi_jVu\|_{L^p(|z|<\rho)}
\\
&\leq
\|V\|_{L^2(B(0,\rho))}
\||z|^{-\tau}\Psi_j u\|_{L^q(|z|<\rho)},
\end{align*}
where we used H\"{o}lder's inequality with
$
\frac1 p=\frac1 2+\frac1 q.
$
Therefore,
\begin{align*}
\||z|^{-\tau}\Psi_j u\|_{L^q(B(0,\rho))}
&\leq
C\||z|^{-\tau}(\bar\partial\Psi_j)f\|_{L^p(\mathbb C)}
\\
&\quad+
C\|V\|_{L^2(B(0,\rho))}
\||z|^{-\tau}\Psi_j u\|_{L^q(B(0,\rho))}+\||z|^{-\tau}\Psi_j\bar\partial f\|_{L^p(|z|>\rho)}.
\end{align*}
By the choice of $\rho$, the last term can be absorbed into the left-hand side, and hence
\[
\||z|^{-\tau}\Psi_j u\|_{L^q(B(0,\rho))}
\leq
C\||z|^{-\tau}(\bar\partial\Psi_j)f\|_{L^p(\mathbb C)}
+
C\||z|^{-\tau}\Psi_j \bar\partial f\|_{L^p(|z|>\rho)} .
\]

The support of $\nabla\Psi_j$ is contained in
$
\left\{z:\frac1j<|z|<\frac2j\right\},
$
and therefore
\begin{align*}
\||z|^{-\tau}(\bar\partial\Psi_j)f\|_{L^p(\mathbb C)}
&\leq
Cj
\||z|^{-\tau}u\|_{L^p(|z|<2/j)}
\\
&\leq
Cj^{\tau+1}
\|u\|_{L^p(|z|<2/j)} .
\end{align*}
Since $u$ is flat at the origin, the above quantity converges to zero as
$j\to\infty$ for every fixed $\tau$.
Passing to the limit $j\to\infty$, we obtain
$$
\||z|^{-\tau}u\|_{L^q(B(0,\rho))}
\leq
C\rho^{-\tau}
\|\bar\partial f\|_{L^p(|z|>\rho)}\leq C\rho^{-\tau}
\| f\|_{W^{1,p}}.
$$
Equivalently,
\[
\left\|
\left(\frac{\rho}{|z|}\right)^\tau u
\right\|_{L^q(B(0,\rho))}
\leq
C\| f\|_{W^{1,p}}.
\]
Since $\rho/|z|>1$ in $B(0,\rho)$, letting $\tau\to\infty$ and applying the
monotone convergence theorem gives
$
u=0$ in $B(0,\rho)$.
This implies $u\equiv0$ by a standard connectedness argument.

\end{proof}

\section{Proof of Theorem~\ref{thm2}}
We now prove the sharpness of the $L^2$ integrability assumption on the
potential. Our construction produces a nontrivial smooth function that
vanishes to infinite order at the origin, together with a potential
belonging to $L^p(\mathbb C)$ for every $0<p<2$, but not to $L^2(\mathbb C)$.
We use the convention $\bar\partial=\frac12(\partial_x+i\partial_y)$,
where $z=x+iy$. For a multi-index
$\beta=(\beta_1,\beta_2)\in\mathbb N_0^2$, write
$|\beta|=\beta_1+\beta_2$ and
$D^\beta=\partial_x^{\beta_1}\partial_y^{\beta_2}$, where $\mathbb N_0=\{0,1,2,\ldots\}$.

Let
\[
K(z)=\frac1{z-1},\qquad
p_{n-1}(z)=-\sum_{k=0}^{n-1}z^k,\qquad n\geq1.
\]
Since $K(z)=-\sum_{k=0}^\infty z^k$ for $|z|<1$, we have
\[
f_n(z)=K(z)-p_{n-1}(z)=\frac{z^n}{z-1}.
\]
Thus $f_n$ is holomorphic near the origin, vanishes there to order $n$,
and has a simple pole at $z=1$. Moreover,
\[
f_n(z)=O(|z|^n)\quad(z\to0),\qquad
f_n(z)=O(|z|^{n-1})\quad(|z|\to\infty).
\]
Fix a cut-off function $\Psi\in C^\infty(\mathbb C)$ satisfying
\[
0\leq\Psi\leq1,\qquad
\Psi(z)=1\quad(|z|\leq\tfrac12),\qquad
\Psi(z)=0\quad(|z|\geq2),
\]
and put $T=\{z\in\mathbb C:\tfrac12<|z|<2\}$.

We first establish the following gluing lemma.
\begin{lem}\label{lem:gluing}
Let $n\geq1$, $r\geq8$, and $0<s\leq\frac 1 8$. There exists
$g_n^{r,s}\in C^\infty(\mathbb C\setminus\{r,s\})$ such that
\[
g_n^{r,s}(z)=r^n f_n(z/r)\quad(|z|\geq2),\qquad
g_n^{r,s}(z)=-s^n f_{n+1}(z/s)\quad(|z|\leq\tfrac12).
\]
Moreover, $g_n^{r,s}$ vanishes only at the origin, and
\begin{equation}\label{eq:gluing-relative}
\operatorname{supp}~(\bar\partial g_n^{r,s})
\subset\overline T,
\qquad
|\bar\partial g_n^{r,s}|
\leq C_0(r^{-1}+s)|g_n^{r,s}|,
\end{equation}
where $C_0>0$ depends only on $\Psi$.
In particular, for any $\epsilon>0$,
the coefficient $C_0(r^{-1}+s)$ can be made
smaller than $\epsilon$ by choosing $r$
sufficiently large and $s$ sufficiently small.
\end{lem}

\begin{proof}
Set
\[
A(z)=r^n f_n(z/r)=-\frac{z^n}{1-z/r},\qquad
B(z)=-s^n f_{n+1}(z/s)=-\frac{z^{n+1}}{z-s},
\]
and let $g_n^{r,s}=(1-\Psi)A+\Psi B$ on
$\mathbb C\setminus\{r,s\}$.  For $z\in T$, write
\[
g_n^{r,s}(z)=-z^nH(z),\]
where
$H(z)=\frac{1-\Psi(z)}{1-z/r}+\frac{\Psi(z)}{1-s/z}.
$
As $|\frac z r|\leq\frac 1 4$ and $|\frac s z|\leq\frac 1 4$ on $T$, we have
\[
\begin{aligned}
|H(z)-1|
&\leq(1-\Psi(z))\frac{|z|/r}{1-|z|/r}
   +\Psi(z)\frac{s/|z|}{1-s/|z|}\\
&\leq\frac13.
\end{aligned}
\]
Consequently,
\begin{equation}\label{eq:lower-transition}
|g_n^{r,s}(z)|\geq\frac23|z|^n,\qquad z\in T.
\end{equation}
Moreover,
\[
|A(z)+z^n|\leq\frac{8}{3r}|z|^n,\qquad
|B(z)+z^n|\leq\frac{8s}{3}|z|^n.
\]
This implies that  $|B-A|\leq\frac83(r^{-1}+s)|z|^n$ on $T$.
Since $A$ and $B$ are holomorphic on $T$, we have
\[
\bar\partial g_n^{r,s}=(\bar\partial\Psi)(B-A).
\]
Combining this identity with \eqref{eq:lower-transition} gives
\[
|\bar\partial g_n^{r,s}(z)|\leq
4\|\bar\partial\Psi\|_\infty
(r^{-1}+s)|g_n^{r,s}(z)|,\qquad z\in T.
\]
Outside $\overline T$, $g_n^{r,s}$ agrees with
the holomorphic functions $B$ and $A$, respectively.
Thus
\[
\operatorname{supp}(\bar\partial g_n^{r,s})
\subset\overline T.
\]
Combining this with the estimate on $T$ yields
\eqref{eq:gluing-relative} with
$C_0=1+4\|\bar\partial\Psi\|_\infty$.
Finally, \eqref{eq:lower-transition} and the
explicit formulas for $A$ and $B$ yield that
the origin is the only zero of $g_n^{r,s}$
in its domain.
\end{proof}

We next regularize the Cauchy kernel near its pole at $z=s$, while preserving control of the relative $\bar\partial$-derivative.
Choose $h\in C^\infty([0,\infty))$ satisfying
\[
h(t)=t^2\quad(0\leq t\leq\tfrac14),
\qquad
h(t)=1\quad(t\geq1),
\qquad h(t)>0\quad(t>0).
\]
Then $C_h=\sup_{t>0}\frac{t|h'(t)|}{h(t)}<\infty.$
For $a>0$, let
\begin{equation*}
Q_a(w)=
\begin{cases}
\displaystyle\frac{h(|w|^2/a^2)}{w},&w\neq0,\\[5pt]
0,&w=0.
\end{cases}
\end{equation*}
Since $Q_a(w)=\frac{w\overline w^{\,2}}{a^4}$ for $|w|\leq \frac a 2$, we have $Q_a\in C^\infty(\mathbb C)$ and
$Q_a(0)=\bar\partial Q_a(0)=0$.
Moreover,
\[
Q_a(w)=\frac1w\quad (|w|\geq a),
\qquad
\bar\partial Q_a(w)
=a^{-2}h'(|w|^2/a^2).
\]
It follows that, for $w\neq0$,
\begin{equation}\label{eq:Q-relative}
\frac{|\bar\partial Q_a(w)|}{|Q_a(w)|}
=\frac1{|w|}\,
\frac{(|w|^2/a^2)|h'(|w|^2/a^2)|}{h(|w|^2/a^2)}
\leq\frac{C_h}{|w|}.
\end{equation}
In particular, the ratio equals $\frac {2} {|w|}$ for
$0<|w|<\frac a 2$.
For every multi-index $\beta$, differentiation of
$Q_a(w)=a^{-1}Q_1(w/a)$ gives
\begin{equation*}
\sup_{|w|\leq a}|D^\beta Q_a(w)|
\leq C_\beta a^{-1-|\beta|},
\end{equation*}
where $C_\beta>0$ is independent of $a$.

We have the following uniform derivative estimates for the regularized function.
\begin{lem}\label{lem:regularized}
Let $n\geq1$, $r\geq8$, $0<s\leq\frac18$, and
$0<\rho<\frac{1}{n+2}$. Set $a=\rho s$ and
$D=B(s,a)$. Define
\begin{equation}\label{eq:G-definition}
G_{n,\rho}^{r,s}(z)=
\begin{cases}
-z^{n+1}Q_a(z-s), & z\in D,\\
g_n^{r,s}(z), & z\notin D,
\end{cases}
\end{equation}
for $z\in\mathbb C\setminus\{r\}$.
Then $G_{n,\rho}^{r,s}\in
C^\infty(\mathbb C\setminus\{r\})$
and agrees with $g_n^{r,s}$ outside $D$.
Its only zeros are $0$ and $s$, with
\[
G_{n,\rho}^{r,s}(s)
=\bar\partial G_{n,\rho}^{r,s}(s)=0.
\]
In addition,
\begin{equation*}
\frac{|\bar\partial G_{n,\rho}^{r,s}(z)|}
{|G_{n,\rho}^{r,s}(z)|}
\leq\frac{C_h}{|z-s|},
\qquad z\in D\setminus\{s\},
\end{equation*}
and the ratio equals $\frac{2}{|z-s|}$
whenever $0<|z-s|<\frac{a}{2}$.

Furthermore, for every multi-index
$\beta\in\mathbb N_0^2$, there exists
a constant $A_{n,\beta,\rho}>0$,
depending only on $n,\beta,\rho$ and
the fixed functions $\Psi,h$, but
independent of $r,s$, such that
\begin{equation}\label{eq:G-derivative}
|D^\beta G_{n,\rho}^{r,s}(z)|
\leq A_{n,\beta,\rho}|z|^{n-|\beta|}
\end{equation}
whenever
$
\frac{n+1}{n+2}s
\leq |z|\leq\frac{n}{n+1}r.
$
\end{lem}

\begin{proof}
Since $a=\rho s<\frac{s}{n+2}$, we have
\[
\overline D\subset
\left\{z:\frac{n+1}{n+2}s<|z|<2s\right\}
\subset\left\{z:|z|<\frac12\right\}.
\]
Thus
$
g_n^{r,s}(z)=-\frac{z^{n+1}}{z-s}
$
in a neighborhood of $\overline D\setminus\{s\}$.
Since $Q_a$ is smooth and agrees with
$\frac 1 w$ for $|w|\geq a$, the two expressions in \eqref{eq:G-definition}
agree in a neighborhood of $\partial D$.
 Hence
$G_{n,\rho}^{r,s}\in
C^\infty(\mathbb C\setminus\{r\})$.
For $|z-s|\leq\frac a2$, we have
\[
G_{n,\rho}^{r,s}(z)
=-\frac{z^{n+1}(z-s)(\overline z-s)^2}{a^4}.
\]
In particular,
\[
G_{n,\rho}^{r,s}(s)
=\bar\partial G_{n,\rho}^{r,s}(s)=0.
\]
Since $0\notin D$ and $Q_a(w)\neq0$
for $w\neq0$, the only zero of
$G_{n,\rho}^{r,s}$ in $D$ is $s$.
Outside $D$, the gluing lemma shows
that its only zero is the origin.

For $z\in D\setminus\{s\}$, the
holomorphicity of $z^{n+1}$ gives
\[
\bar\partial G_{n,\rho}^{r,s}(z)
=-z^{n+1}\bar\partial Q_a(z-s).
\]
Consequently, by \eqref{eq:Q-relative},
\[
\frac{|\bar\partial G_{n,\rho}^{r,s}(z)|}
{|G_{n,\rho}^{r,s}(z)|}
=
\frac{|\bar\partial Q_a(z-s)|}
{|Q_a(z-s)|}
\leq\frac{C_h}{|z-s|}.
\]
Moreover, the ratio equals
$\frac{2}{|z-s|}$ whenever
$0<|z-s|<\frac a2$.

It remains to prove \eqref{eq:G-derivative}.
We consider four cases.
First, if
\[
\frac{n+1}{n+2}s\leq|z|\leq2s,
\]
writing $z=s\zeta$ yields
\[
G_{n,\rho}^{r,s}(s\zeta)
=-s^n\zeta^{n+1}Q_\rho(\zeta-1).
\]
Since $Q_\rho$ is smooth and
$\frac{n+1}{n+2}\leq|\zeta|\leq2$,
differentiation gives
\[
|D^\beta G_{n,\rho}^{r,s}(z)|
\leq C_{n,\beta,\rho}s^{n-|\beta|}
\leq C'_{n,\beta,\rho}|z|^{n-|\beta|}.
\]

Next, for $2s\leq|z|\leq\frac12$,
we have
\[
G_{n,\rho}^{r,s}(z)
=-\frac{z^{n+1}}{z-s}.
\]
Since $|z-s|\geq\frac{|z|}{2}$,
Leibniz's rule gives
\[
|D^\beta G_{n,\rho}^{r,s}(z)|
\leq C_{n,\beta}|z|^{n-|\beta|}.
\]

For $\frac12\leq|z|\leq2$, write
\[
G_{n,\rho}^{r,s}(z)=-z^nH(z),
\]
where
\[
H(z)=\frac{1-\Psi(z)}{1-z/r}
+\frac{\Psi(z)}{1-s/z}.
\]
Since $r\geq8$ and $s\leq\frac18$,
all derivatives of $H$ are uniformly
bounded on this annulus, independently
of $r$ and $s$. Thus
\[
|D^\beta G_{n,\rho}^{r,s}(z)|
\leq C_{n,\beta}|z|^{n-|\beta|}.
\]

Finally, for
$2\leq|z|\leq\frac{n}{n+1}r$,
we have
\[
G_{n,\rho}^{r,s}(z)
=-\frac{z^n}{1-z/r}.
\]
Since
\[
\left|1-\frac zr\right|
\geq\frac{1}{n+1},
\qquad
\frac1r\leq\frac1{|z|},
\]
repeated differentiation implies
\[
|D^\beta G_{n,\rho}^{r,s}(z)|
\leq C_{n,\beta}|z|^{n-|\beta|}.
\]
Combining these estimates proves
\eqref{eq:G-derivative}, with a constant
independent of $r$ and $s$.
\end{proof}


We now combine the local constructions at decreasing scales to
prove Theorem~\ref{thm2}. The matching formulas ensure smoothness
across the annular interfaces, while the uniform derivative bounds
and relative $\bar\partial$-estimates provide flatness and
integrability, respectively.

\begin{proof}[Proof of Theorem~\ref{thm2}]
Fix $0<b<\frac{1}{8}$ and set
$\rho_j=\frac{b}{j+1}$ for $j\geq1$.
Let $r_1=1$, $P_0=1$, and
$P_n=\prod_{j=1}^{n}r_j$ for $n\geq1$.
Starting from $n=2$, we inductively choose
positive numbers $\eta_n$ and $r_n$ such that
\begin{equation}\label{eq:separation}
0<r_n<\frac{\eta_n}{8},\qquad
0<\eta_n<\frac{r_{n-1}}{8},
\end{equation}
\begin{equation}\label{eq:normalized-error}
C_0\left(
\frac{\eta_n}{r_{n-1}}+\frac{r_n}{\eta_n}
\right)
\leq 2^{-(n-1)},
\end{equation}
and, for every $\beta\in\mathbb N_0^2$
with $|\beta|<\frac{n}{2}$,
\begin{equation}\label{eq:flatness-choice}
\frac{A_{n,\beta,\rho_{n+1}}}{P_{n-1}}
r_n^{n-|\beta|-1}\leq 2^{-n}.
\end{equation}
Here $A_{n,\beta,\rho}$ is the constant in
\eqref{eq:G-derivative}, independent of $r$ and $s$.
Such a choice is possible by taking $\eta_n$ and
then $r_n$ sufficiently small, since
$n-|\beta|-1>0$ for $|\beta|<\frac{n}{2}$
and only finitely many such multi-indices occur.
Moreover, \eqref{eq:separation} implies
\[
r_n<\frac{r_{n-1}}{64},\qquad n\geq2,
\]
so $r_n\to0$ geometrically.

For $n\geq1$, define
\begin{equation*}
E_{n+1}=\left\{z\in\mathbb C:
\frac{n+1}{n+2}r_{n+1}\leq|z|\leq
\frac{n}{n+1}r_n\right\}.
\end{equation*}
These closed annuli cover $0<|z|\leq\frac12$
and intersect only along their common boundary circles.
Define
\begin{equation}\label{eq:outer-u}
u(z)=-zQ_{\rho_1}(z-1),
\qquad |z|\geq\frac12,
\end{equation}
and, for $z\in E_{n+1}$, set
\begin{equation*}
u(z)=c_{n+1}
G_{n,\rho_{n+1}}^{
r_n/\eta_{n+1},\,r_{n+1}/\eta_{n+1}}
\left(\frac{z}{\eta_{n+1}}\right),
\qquad
c_{n+1}=\frac{(-1)^n\eta_{n+1}^{n}}{P_n}.
\end{equation*}
Since
\[
\frac{r_n}{\eta_{n+1}}>8,\qquad
\frac{r_{n+1}}{\eta_{n+1}}<\frac18,
\qquad
\rho_{n+1}<\frac{1}{n+2},
\]
both lemmas apply. The remaining pole of the
regularized function lies at $z=r_n$,
outside $E_{n+1}$, so $u$ is well defined
on each annulus.

We verify that the definitions agree across the annular interfaces. Near $|z|=\frac{n}{n+1}r_n$,
the function on $E_{n+1}$ agrees with the outer holomorphic model. Thus
\begin{equation}\label{eq:matching}
\begin{aligned}
u(z)
&=c_{n+1}\left(\frac{r_n}{\eta_{n+1}}\right)^n
f_n(z/r_n)\\
&=\frac{(-1)^n r_n^{n-1}}{P_{n-1}}
f_n(z/r_n).
\end{aligned}
\end{equation}
For $n\geq2$, this agrees with the inner
holomorphic model on $E_n$, since
\[
-c_n\left(\frac{r_n}{\eta_n}\right)^{n-1}
=\frac{(-1)^n r_n^{n-1}}{P_{n-1}}.
\]
Moreover, the regularization disc centered
at $r_n$ does not meet the interface, as
$\rho_n<\frac{1}{n+1}$.
For $n=1$, \eqref{eq:matching} becomes
$-f_1(z)$, agreeing with \eqref{eq:outer-u}
near $|z|=\frac12$.
Hence the definitions agree in neighborhoods
of all interfaces, and
$
u\in C^\infty(\mathbb C\setminus\{0\}).
$

We now extend $u$ smoothly across the origin.
By the chain rule and \eqref{eq:G-derivative},
for $z\in E_{n+1}$,
\begin{equation*}
|D^\beta u(z)|
\leq \frac{A_{n,\beta,\rho_{n+1}}}{P_n}
|z|^{n-|\beta|}.
\end{equation*}
For $n\geq2$ and $|\beta|<\frac{n}{2}$,
since $|z|\leq r_n$ and $P_n=P_{n-1}r_n$,
\eqref{eq:flatness-choice} yields
\[
|D^\beta u(z)|
\leq
\frac{A_{n,\beta,\rho_{n+1}}}{P_{n-1}}
r_n^{n-|\beta|-1}
\leq 2^{-n}.
\]
Consequently, for every multi-index $\beta$,
$
\lim_{z\to0}D^\beta u(z)=0.
$
Set $u(0)=0$ and extend each $D^\beta u$
continuously to the origin by zero,
denoting the extension by $F_\beta$.
By the fundamental theorem of calculus,
\[
F_\beta(z)
=\int_0^1 \nabla D^\beta u(tz)\cdot z\,dt,
\]
and hence
\[
\frac{|F_\beta(z)-F_\beta(0)|}{|z|}
\leq
\sup_{0<|w|\leq|z|}
|\nabla D^\beta u(w)|
\longrightarrow0
\quad\text{as }z\to0.
\]
It follows inductively that
$u\in C^\infty(\mathbb C)$ and
$D^\beta u(0)=0$ for every $\beta$.
Thus $u$ vanishes to infinite order at
the origin.

We now establish the required estimate for the potential.
For $n\geq1$, let
\[
B_{n+1}=\left\{z\in\mathbb C:
\frac{\eta_{n+1}}{2}<|z|<2\eta_{n+1}\right\},
\qquad
D_n=B(r_n,\rho_n r_n).
\]
By the separation conditions, these annuli and discs
are pairwise disjoint and contained in $B(0,2)$.
Moreover, $u$ is holomorphic outside their union
and the origin, and the two lemmas imply that
its zero set is precisely
$\{0\}\cup\{r_n:n\geq1\}$.
Define
\begin{equation*}
V(z)=
\begin{cases}
\frac{|\bar\partial u(z)|}{|u(z)|},
&u(z)\neq0,\\[5pt]
0,&u(z)=0.
\end{cases}
\end{equation*}
Then $V$ is nonnegative and measurable, and
\[
|\bar\partial u(z)|=V(z)|u(z)|,
\qquad z\in\mathbb C.
\]

On $B_{n+1}$, the rescaled gluing estimate
and \eqref{eq:normalized-error} yield
\[
V(z)\leq\frac{C_0}{\eta_{n+1}}
\left(\frac{\eta_{n+1}}{r_n}
+\frac{r_{n+1}}{\eta_{n+1}}\right)
\leq\frac{2^{-n}}{\eta_{n+1}}.
\]
Since $|B_{n+1}|\leq4\pi\eta_{n+1}^2$,
for every $0<p<2$,
\begin{equation*}
\sum_{n=1}^\infty
\int_{B_{n+1}}V(z)^p\,dv_z
\leq4\pi\sum_{n=1}^\infty
2^{-np}\eta_{n+1}^{2-p}<\infty.
\end{equation*}
On $D_n\setminus\{r_n\}$, the regularization
estimate gives
\begin{equation}\label{eq:D-estimate}
V(z)\leq\frac{C_h}{|z-r_n|},
\end{equation}
where the case $n=1$ follows directly
from \eqref{eq:outer-u}. Consequently,
\begin{equation*}
\int_{D_n}V(z)^p\,dv_z
\leq\frac{2\pi C_h^p}{2-p}
(\rho_n r_n)^{2-p}.
\end{equation*}
Since $r_n\leq64^{1-n}$ and $\rho_n<1$,
\[
\sum_{n=1}^\infty
(\rho_n r_n)^{2-p}<\infty.
\]
As $V$ vanishes outside these regions,
we conclude that
\[
V\in\bigcap_{0<p<2}L^p(\mathbb C).
\]

On the other hand, the exact formula for
$Q_a$ yields
\[
V(z)=\frac{2}{|z-r_n|}
\]
whenever
$
0<|z-r_n|<\frac{\rho_n r_n}{2}.
$
In particular,
\[
\int_{D_1}V(z)^2\,dv_z
\geq8\pi\int_0^{\rho_1/2}
\frac{dt}{t}=\infty.
\]
Thus $V\notin L^2(\mathbb C)$.

For general $N\geq1$, set
$U(z)=(u(z),0,\ldots,0)\in\mathbb C^N$.
Then $U$ is smooth, nontrivial, and vanishes to infinite order at
$0$. Moreover,
\[
|\bar\partial U(z)|=|\bar\partial u(z)|
=V(z)|u(z)|=V(z)|U(z)|,
\qquad z\in\mathbb C.
\]
This proves the theorem.
\end{proof}

\vspace{1em}
\noindent \textbf{Acknowledgements} The authors would like to thank the reviewers for their careful review of this manuscript and for their insightful comments and constructive suggestions.

\section*{Declarations}

\noindent \textbf{Author Contributions} All authors contributed equally to the writing of this paper. All authors read and approved the final manuscript.

\vspace{1em}

\noindent \textbf{Funding} X. Tang was supported by the Zhejiang Provincial Natural Science Foundation of China (Grant No. LZ24A010004).
 S. He was supported by the National Natural Science Foundation of China (Grant No. 12301115).

\vspace{1em}

\noindent \textbf{Data Availability} This is no data availability.

\vspace{1em}

\noindent \textbf{Conflict of interest} There is no Conflict of interest.


\begin{thebibliography}{s15}
\bibitem{BL}S. Bell, L. Lempert, A  $C^\infty$ Schwarz reflection principle in one and several complex variables.
\emph{J. Differential Geom. }   \textbf{32} (1990), no. 3, 899--915.


\bibitem{CS}
S. Chanillo, E. Sawyer,
Unique continuation for $\Delta+V$ and the C. Fefferman-Phong class.
\emph{Trans. Amer. Math. Soc.} \textbf{318} (1990), 275--300.



\bibitem{JK1985}D. Jerison, C. E. Kenig,
Unique continuation and absence of positive eigenvalues for Schr\"{o}dinger operators.
\emph{Ann. of Math. (2)} \textbf{121} (1985), no. 3, 463--494.

\bibitem{kenig}C. E. Kenig,  J.-N. Wang, Quantitative uniqueness estimates for second order elliptic equations with unbounded drift. \emph{Math. Res. Lett.} \textbf{22} (2015), no. 4, 1159--1175.

\bibitem{MAND}N. Mandache,
A counterexample to unique continuation in dimension two.
\emph{Comm. Anal. Geom.} \textbf{10} (2002), no. 1, 1--10.






\bibitem{PZ}Y. Pan, Y.  Zhang,
Unique continuation for $\bar{\partial}$ with square-integrable potentials.
\emph{New York J. Math.}  \textbf{29} (2023), 402--416.


\bibitem{PZ2}Y. Pan, Y.  Zhang, Unique continuation of Schr\"{o}dinger-type equations for  $\bar{\partial}$.
\emph{Comm. Anal. Geom.} \textbf{33} (2025), no. 8, 1975--2014.

\bibitem{seo}I. Seo,
A remark on unique continuation for the Cauchy-Riemann operator.
\emph{Bull. Korean Math. Soc.}\textbf{ 52} (2015), no. 5, 1753--1757.

\bibitem{S1990}C. Sogge, Christopher,
Strong uniqueness theorems for second order elliptic differential equations.
\emph{Amer. J. Math. }\textbf{112} (1990), no. 6, 943--984.



\bibitem{wolff1990}T. Wolff,
Unique continuation for $|\Delta u|\leq V|\nabla u|$ and related problems.
\emph{Rev. Mat. Iberoamericana} \textbf{6} (1990), no. 3--4, 155--200.

\bibitem{wolff1992}T. Wolff,
Note on counterexamples in strong unique continuation problems.
\emph{Proc. Amer. Math. Soc.} \textbf{114} (1992), no. 2, 351--356.

\end{thebibliography}
\end{document}